\documentclass[12pt]{article}
\usepackage{amsfonts,amssymb}
\usepackage{amsthm}
\newtheorem{theorem}{Theorem}
\newtheorem{lemma}{Lemma}

\newtheorem{corollary}{Corollary}

\newtheorem*{theorem*}{Theorem *}
\title{Any Sub-Riemannian Metric has Points of Smoothness}
\author{A.~A.~Agrachev\thanks{SISSA, Trieste \& Steklov Math. Inst.,
Moscow}}
\date{}
\begin{document}
\maketitle
\begin{abstract} We prove the result stated in the title; it is
equivalent to the existence of a regular point of the
sub-Riemannian exponential mapping. We also prove that the metric
is analytic on an open everywhere dense subset in the case of a
complete real-analytic sub-Riemannian manifold.
\end{abstract}

\subsection*{1. Preliminaries}
Let $M$ be a smooth (i.\,e. $C^\infty$) Riemannian manifold and
$\Delta\subset TM$ a smooth vector distribution on $M$ (a vector
subbundle of $TM$). We denote by $\bar\Delta$ the space of smooth
sections of $\Delta$ that is a subspace of the space
$\mathrm{Vec}M$ of smooth vector fields on $M$. The Lie bracket of
vector fields $X,Y$ is denoted by $[X,Y]$. We assume that $\Delta$
is bracket generating; in other words, $\forall q\in M$,
$$
span\{[X_1,[\cdots,[X_{m-1},X_m]\cdots](q) : X_i\in\bar\Delta,\
i=1,\ldots m,\ m\in\mathbb N\}=T_qM.
$$

Given $q_0,q_1\in M$, we define the space of starting from $q_0$
{\it admissible paths}:
$$
\Omega_{q_0}=\{\gamma\in H^1([0,1],M) : \gamma(0)=q_0,\
\dot\gamma(t)\in\Delta_{\gamma(t)}\ \mathrm{for\ almost\ all}\ t\}
$$
and the sub-Riemannian distance:
$$
\delta(q_0,q_1)=\inf\{\ell(\gamma) : \gamma\in\Omega_{q_0},\
\gamma(1)=q_1\},
$$
where $\ell(\gamma)=\int_0^1|\dot\gamma(t)|\,dt$ is the length of
$\gamma$ and $\Delta_q=\Delta\cap T_qM$. Classical
Rashevskij--Chow theorem implies that $\delta$ is a well-defined
continuous function on $M\times M$.

An admissible path $\gamma$ is a {\it length-minimizer} if
$\ell(\gamma)=\delta(\gamma(0),\gamma(1))$. Given $s\in[0,1]$, we
define a re-scaled path $\gamma^s:t\mapsto\gamma(st),\ t\in[0,1]$.
The re-scaled paths of any length-minimizer are also
length-minimizers. According to the standard Filippov existence
theorem, any $q\in M$ belongs to the interior of the set of points
connected with $q$ by a length-minimizer. If $M$ is a complete
Riemannian manifold, then length-mi\-ni\-mi\-zers connect $q$ with
all points of $M$.

From now on, the point $q_0\in M$ is supposed to be fixed. Note
that $\Omega_{q_0}$ is a smooth Hilbert submanifold of
$H^1([0,1],M)$. A smooth {\it endpoint mapping} $f:\Omega_{q_0}\to
M$ is defined by the formula: $f(\gamma)=\gamma(1)$. Let $J:\Omega
_{q_0}\to\mathbb R$ be the action functional, $J(\gamma)=\frac
12\int_0^1|\dot\gamma(t)|^2\,dt$. The Cauchy--Schwartz inequality
implies that an admissible curve realizes $\min\limits_{\gamma\in
f^{-1}(q)}J(\gamma)$ if and only if this curve is a connecting
$q_0$ with $q$ and parameterized proportionally to the length
length-minimizer. Then, according to the Lagrange multipliers
rule, any starting at $q_0$ and parameterized proportionally to
the length length-minimizer is either a critical point of $f$ or a
solution of the equation
$$
\lambda D_\gamma f=d_\gamma J \eqno (1)
$$
for some $\lambda\in T^*_{\gamma(1)}M$, where $\lambda D_\gamma
f\in T^*_\gamma\Omega _{q_0}$ is the composition of linear
mappings $D_\gamma f:T_\gamma\Omega _{q_0}\to T_{\gamma(1)}M$ and
$\lambda:T_{\gamma(1)}M\to\mathbb R$.

Now let $a\in C^1(M)$ and a curve $\gamma$ realizes
$\min\limits_{\gamma\in\Omega
_{q_0}(M,\Delta)}\left(J(\gamma)-a(\gamma(1))\right)$; then
$D_\gamma J-d_{\gamma(1)}aD_\gamma f=0$. Hence $\gamma$ satisfies
(1) with $\lambda=d_{\gamma(1)}a$.

Solutions of (1) are called {\it normal (sub-Riemannian)
geodesics} while critical points of $f$ are {\it abnormal
geodesics}. If a geodesic satisfies (1) for at least two different
$\lambda$, then it is simultaneously normal and abnormal; all
other geodesics are either strictly normal or strictly abnormal.

The {\it sub-Riemannian Hamiltonian} is a function
$h:T^*M\to\mathbb R$ defined by the formula: $
h(\lambda)=\max\limits_{v\in\Delta_q}\left(\langle\lambda,v\rangle-\frac
12|v|^2\right),\ \lambda\in T^*_qM,\ q\in M. $ We denote by $\vec
h$ the associated to $h$ Hamiltonian vector field on $T^*M$. A
pair $(\gamma,\lambda)$ satisfies (1) if and only if there exists
a solution $\psi$ of the Hamiltonian system $\dot\psi=\vec
h(\psi)$ such that $\psi(1)=\lambda$ and $\psi(t)\in
T^*_{\gamma(t)}M,\ \forall t\in[0,1]$; this fact is a very special
case of the Pontryagin maximum principle.

Note that $h\bigr|_{T^*_qM}$ is a nonnegative quadratic form on
$T^*_qM$ whose kernel equals $\Delta^\perp_q=\{\lambda\in
T^*_qM:\lambda\perp\Delta_q\}$. Given $\xi\in T^*_{q_0}M$ we
denote by $\gamma_\xi$ normal geodesic that is the projection to
$M$ of the solution of the Cauchy problem: $\dot\psi=\vec
h(\psi),\ \psi(0)=\xi$. This notation is well-coordinated with
re-scalings: $\gamma^s_\xi=\gamma_{s\xi},\ \forall s\in[0,1]$.
Finally, we define the {\it exponential map} $\mathcal
E:\xi\mapsto\gamma_\xi(1)$; this is a smooth map of a neighborhood
of $\Delta^\perp_{q_0}$ in $T^*_{q_0}M$ to $M$ and $\mathcal
E(\Delta^\perp_{q_0})=q_0$.

\subsection*{2. Statements} A point $q\in M$ is called a {\it
smooth point} (for the triple $(M,\Delta,q_0)$) if
$\exists\,\xi\in T_{q_0}^*M$ such that $q=\mathcal E(\xi)$, $\xi$
is a regular point of $\mathcal E$ and $\gamma_\xi$ is a unique
length-minimizer connecting $q_0$ with $q$. We denote by $\Sigma$
the set of all smooth points and assume (for all statements of
this section) that $M$ is a complete Riemannian manifold.

\begin{theorem} $\Sigma$ is an open everywhere dense subset of $M$.
\end{theorem}

The term ``smooth point" is justified by the following fact, which
should be well-known to the experts even if it is not easy to find
an appropriate reference. A very close statement is contained in
\cite{Ja}.

\begin{theorem*}
\begin{description} \item[\rm(i)] If $q\in\Sigma$,
then the sub-Riemannian distance $\delta$ is smooth in a
neighborhood of $(q_0,q)$. If, additionally, $M$ and $\Delta$ are
real-analytic, then $\delta$ is analytic in a neighborhood of
$(q_0,q)$.
\item[\rm(ii)] If $\delta$ is $C^2$ in a neighborhood
of $(q_0,q)$, then $q\in\Sigma$ and $\delta$ is actually smooth at
$(q_0,q)$.
\end{description}
\end{theorem*}

\begin{corollary} Sub-Riemannian distance $\delta$ is smooth
on an open everywhere dense subset $\mathcal S\subset M\times M$.
Moreover, $\mathcal S \cap\{q_0\}\times M$ is everywhere dense in
$\{q_0\}\times M,\ \forall q_0\in M$. If $M$ and $\Delta$ are
real-analytic, then $\delta$ is analytic on $\mathcal S$.
\end{corollary}

\subsection*{2. Proofs}
We have: $\mathrm{im}D_\xi\mathcal
E\subset\mathrm{im}D_{\gamma_\xi}f,\ \forall\xi\in T^*_{q_0}M$,
since $\mathcal E(\xi)\equiv f(\gamma_\xi)$. Given a normal
geodesic $\gamma$ we say that the point $\gamma(1)$ is conjugate
to $q_0$ along $\gamma$ if $\mathrm{im}D_\xi\mathcal
E\ne\mathrm{im}D_{\gamma_\xi}f$ for some $\xi$ such that
$\gamma=\gamma_\xi$. The following three properties of conjugate
points are specifications of more general facts whose proofs can
be found in \cite[Ch.\,21]{AgSa}.

i) If $\mathcal E(s\xi)=\gamma_\xi(s)$ is not conjugate to $q_0,\
\forall s\in [0,1]$, then $\gamma_\xi$ is strictly shorter than
any other connecting $q_0$ with $\gamma_\xi(1)$ admissible path
from a $C^0$-neighborhood of $\gamma_\xi$.

ii) If $\xi$ is a regular point of $\mathcal E$ and $\gamma_\xi$
is strictly shorter than any other connecting $q_0$ with
$\gamma_\xi(1)$ admissible path from a $C^\infty$-neighborhood of
$\gamma_\xi$, then $\gamma_\xi(s)$ is not conjugate to $q_0$ along
$\gamma_{s\xi}$, $\forall s\in [0,1]$.

iii) The set $\{s\in[0,1]:\gamma_\xi(s)\ \mathrm{is\ conjugate\
to}\ q_0\ \mathrm{along}\ \gamma_{s\xi}\}$ is a closed subset of
$[0,1]$ which does not contain 0. Moreover, this is a finite
subset of $[0,1]$ if $M$ and $\Delta$ are real-analytic.

We say that that $q\in M$ is a {\it RT-point} (after Rifford and
Trelat) if $q=\mathcal E(\xi)$ for some $\xi\in T^*M$ such that
$\gamma_\xi$ is a unique length-minimizer connecting $q_0$ with
$q$. Obviously, any smooth point is a RT-point but not vice versa!
In particular, a normal geodesic $\gamma_\xi$ from the definition
of the RT-point can be also abnormal. If $M$ is complete, then the
set of RT-points is everywhere dense in $M$. This fact is proved
in \cite{RiTr}; the proof is simple and we present it here.

Given an open subset $O\subset M$ we denote by $RT_O$ the set of
all RT-points of $O$. We have to show that $RT_O$ is not empty.
Let $a:O\to\mathbb R$ be a smooth function such that
$a^{-1}([c,\infty))$ is compact for any $c\in\mathbb R$. Then the
function $q'\mapsto\delta(q_0,q')-a(q'),\ q'\in O,$ attains
minimum at some point $q\in O$. Hence any connecting $q_0$ with
$q$ length-minimizer $\gamma$ satisfies equation (1) with
$\lambda=d_qa$. Then $\gamma$ is the projection to $M$ of the
solution to the Cauchy problem $\dot\psi=\vec h(\psi),\
\psi(1)=d_qa$; in other words, $\gamma=\gamma_{\psi(0)}$.

Now we prove that $\Sigma$ is everywhere dense in $M$.
Suppose that there exists an open
subset $O\subset M$ such that any point of $RT_O$ is connected
with $q_0$ by an abnormal length-minimizer. Given $q\in RT_O$ we
set $\mathrm{rk}(q)=\mathrm{rank}\,D_\gamma f$, where $\gamma$ is
the length-mi\-ni\-mi\-zer connecting $q_0$ with $q$. Finally, let
$k_O=\max\limits_{q\in RT_O}\mathrm{rk}(q)$. According to our
assumption, $k_O<\dim M$.

Now take $\hat q\in RT_O$ such that $\mathrm{rk}(\hat q)=k_O$.
Then $\mathrm{rk}(q)=k_O$ for any sufficiently close to $\hat q$
point $q\in RT_O$. Indeed, take a convergent to $\hat q$ sequence
$q_n\in RT_O,\ n=1,2,\ldots$ . Let $\gamma_n$ be the
length-mi\-ni\-mi\-zer connecting $q_0$ with $q_n$ and
$\hat\gamma$ be the unique length-mi\-ni\-mi\-zer connecting $q$
with $\hat q$. The uniqueness property and compactness of the
space of length-mi\-ni\-mi\-zers (see \cite{com}) imply that
$\gamma_n\to\hat\gamma$ in the $H^1$-topology and
$D_{\gamma_n}f\to D_{\hat\gamma}f$ as $n\to\infty$. Hence
$\mathrm{rk}(\hat q)\le\mathrm{rk}(q_n)$ for all sufficiently big
$n$.

Now, if necessary, we can substitute $O$ by a smaller open subset
and assume, without lack of generality, that $\mathrm{rk}(q)=k_O,\
\forall q\in O$. Given $q\in RT_O$ and the connecting $q_0$ with
$q$ length-mi\-ni\-mi\-zer $\gamma$ we set
$$
\Pi_q=\{\xi\in T^*_{q_0}M :\gamma_\xi=\gamma\}.
$$
It is easy to see that $\Pi_q$ is an affine subspace of
$T^*_{q_0}M$; moreover, $\xi\in\Pi_q$ if and only if
$\lambda=e^{\vec h}(\xi)$ satisfies (1), where $e^{t\vec
h}:T^*M\to T^*M,\ t\in\mathbb R$, is the generated by $\vec h$
Hamiltonian flow. The already used compactness--uniqueness
argument implies that the affine subspace $\Pi_q\subset
T^*_{q_0}M$ continuously depends on $q\in RT_O$.

Consider again $\hat q=\hat\gamma(1)\in RT_O$ and take a
containing $\hat\xi$ and transversal to $\Pi_{\hat q}$ $(\dim
M-k_O)$-dimensional ball $\mathfrak B$ in $M$. There exists a
neighborhood $\hat O$ of $\hat q$ such that $\Pi_q\cap \mathfrak
B\ne\emptyset,\ \forall q\in\hat O\cap RT_O$. Hence any
sufficiently close to $\hat q$ element of $RT_O$ belongs to the
compact zero measure subset $\mathcal E(\mathfrak B)$. We obtain a
contradiction with the fact that $RT_O$ is everywhere dense in
$M$.

 This contradiction proves that the set
$$RT'_O\stackrel{def}{=}\{q\in RT_O: \mathrm{rk}(q)=\dim M\}$$ of
RT-points connected with $q_0$ by a strictly normal
length-minimizer is everywhere dense in $M$. Now we are going to
find a smooth point arbitrary close to the given point $q\in
RT_O'$. According to the definitions, a point $q'\in RT_O'$ is
smooth if and only if it is not conjugate to $q_0$ along a
length-minimizer $\gamma_\xi$, where $q'=\gamma_\xi(1)$. If $M$
and $\Delta$ are real-analytic, then $\gamma(s)$ is a smooth point
for any sufficiently close to 1 number $s<1$ (see iii)). In the
general $C^\infty$ situation we need the following

\begin{lemma}[L. Rifford] Let $q\in RT'_O$;
then the function
$$
q'\mapsto\delta(q_0,q') \eqno (2)
$$
is Lipschitz in a neighborhood of $q$.
\end{lemma}
{\bf Proof.} In order to prove this local statement, we may fix
some local coordinates in $\Omega_{q_0}$ and $M$ and make all the
computations under assumption that $\Omega_{q_0}$ is a Hilbert
space and $M=\mathbb R^n$. Let $\gamma_\xi$ be the connecting
$q_0$ with $q$ length-minimizer; then $\gamma_\xi$ is a regular
point of the endpoint map $f:\gamma\mapsto\gamma(1),\
\gamma\in\Omega_{q_0}$. Let $B_\gamma(\varepsilon)$ and
$B_q(\varepsilon)$ be the centered at $\gamma$ and $q$ radius
$\varepsilon$ balls in $\Omega_{q_0}$ and $\mathbb R^n$. The
implicit function theorem implies the existence of constants
$c,\alpha>0$ and a neighborhood $\mathcal O_{\gamma_\xi}$ of
$\gamma_\xi$ in $\Omega_{q_0}$ such that
$B_{\gamma(1)}(\varepsilon)\subset
f\left(B_\gamma(c\varepsilon)\right)$ for any $\gamma\in\mathcal
O_{\gamma_\xi}$ and any $\varepsilon\in(0,\alpha]$. Then $
\delta(q_0,q')\le\ell(\gamma)+\bar c|\gamma(1)-q'|,$ for some
constant $\bar c>0$ and any $\gamma\in\mathcal O_{\gamma_\xi},\
q'\in B_{\gamma(1)}(\alpha)$.

Now take a neighborhood $U_q$ of $q$ in $M$ such that $U_q\subset
B_q(\alpha),\ RT_{U_q}=RT'_{U_q}$ and the connecting $q_0$ with
points from $RT'_{U_q}$ length-minimizers belong to $\mathcal
O_{\gamma_\xi}$. Then
$$
\delta(q_0,q_1)-\delta(q_0,q_2)\le\bar c|q_1-q_2|,\quad \forall
q_1,q_2\in RT_{U_q}.
$$
Hence function (2) is Lipschitz on $RT_{U_q}$ and, by the
continuity, on \linebreak $U_q=\overline{RT}_{U_q}. \quad \square$

A Lipschitz function is differentiable almost everywhere. It is
easy to see that any point of differentiability of function (2) is
an RT-point. Indeed, set $\phi(q)=\frac 12\delta^2(q_0,q),\ q\in
M$. The functional $J(\gamma)-\phi(\gamma(1)),\
\gamma\in\Omega_{q_0},$ attains minimum exactly on the
length-minimizers. Hence a connecting $q_0$ with a
differentiability point of $\phi$ length-minimizer $\gamma$ must
satisfy the equation $d_{\gamma(1)}\phi D_\gamma f=d_\gamma J$ and
is thus unique and normal. All RT-points are values of the
exponential map; according to the Sard Lemma, almost all of them
must be regular values and thus smooth points!

\medskip
In the next lemma, we allow to perturb the sub-Riemannian
structure, i.\,e. the Riemannian structure on the given smooth
manifold and the vector distribution $\Delta$. The space of
sub-Riemannian structures (shortly SR-structures) is endowed with
the standard $C^\infty$ topology.

\begin{lemma} Assume that $M$ is complete and $q$ is a smooth
point. Then any sufficiently close to $q$ point is smooth.
Moreover, all sufficiently close to $q$ points remain to be smooth
after a small perturbation of $q_0$ and the SR-structure; the
connecting $q_0$ with $q$ length-minimizer smoothly depends on the
triple $(q,\,q_0,$\,SR-structure).
\end{lemma}
{\bf Proof.} Let $\gamma_\xi$ be the connecting $q_0$ with $q$
length-minimizer. The fact that $\xi$ is a regular point of
$\mathcal E$ allows to find normal geodesics connecting any close
to $q_0$ point with any close to $q$ point for any sufficiently
close to $(M,\Delta)$ SR-structure in such a way that the geodesic
smoothly depends on all the data. It remains to show that the
found geodesic is a unique length-mi\-ni\-mi\-zer connecting
corresponding points!

It follows from property ii) of the conjugate points that
$\gamma_\xi(s)$ is not conjugate to $q_0$ along $\gamma_{s\xi},\
\forall s\in[0,1]$. This fact implies the existence of a
containing $\xi$ Lagrange submanifold $\mathcal L\subset T^*M$
such that $\pi\circ e^{t\vec h}\bigr|_{\mathcal L}$ is a
diffeomorphism of $\mathcal L$ on a neighborhood of
$\gamma_\xi(t),\ \forall t\in [0,1]$, where $\pi:T^*M\to M$ is the
standard projection (see \cite{AgStZe} or \cite[Ch.21]{AgSa}).
Then $\gamma_\xi$ is strictly shorter than any connecting $q_0$
with $q$ admissible path $\gamma$ such that
$$
\gamma(t)\in \pi\circ e^{t\vec h}(\mathcal L),\quad \forall t\in
[0,1]. \eqno (3)
$$
Moreover, compactness of the space of length-mi\-ni\-mi\-zers
implies that $\exists\,\varepsilon>0$ such that
$\ell(\gamma)-\ell(\gamma_\xi)\ge\varepsilon$ for any connecting
$q_0$ with $q$ admissible $\gamma$ which does not satisfy (3). It
remains to mention that construction of $\mathcal L$ survives
small perturbations of the SR-structure and of the initial data
for normal geodesics. Hence found geodesics are indeed unique
length minimizers connecting their endpoints.$\quad \square$

Lemma~2 implies that $\Sigma$ is open as soon as $M$ is complete.
This finishes proof of Theorem~1. Moreover, statement (i) of
Theorem * is also a direct corollary of Lemma~2. Here is the proof
of statement (ii): Let $\phi(q')=\frac 12\delta^2(q_0,q'),\ q'\in
M$, then the functional $\gamma\mapsto J(\gamma)-\phi(\gamma(1)),\
\gamma\in\Omega_{q_0}(M,\Delta)$ attains minimum on the
length-minimizers. Hence $q$ is an RT-point and the connecting
$q_0$ with $q$ length-minimizer is $\gamma_\xi$, where
$\xi=e^{-t\vec h}(d_q\phi)$. Now the mapping $q'\mapsto e^{-t\vec
h}(d_{q'}\phi)$ defines a local inverse of $\mathcal E$ on a
neighborhood of $q=\mathcal E(\xi)$. Hence $\xi$ is a regular
point of $\mathcal E$ and $q\in\Sigma$.

\medskip\noindent {\sl Acknowledgments.} I am grateful to professors
Morimoto, Rifford and Zelenko for very useful discussions.

\medskip\noindent {\sl Remark.} The first version of this paper
with a little bit weaker result to appear in the ``Russian Math.
Dokl.". The improvement made in the present updated version
concerns the density of the set of smooth points in the general
$C^\infty$ case. The density was originally proved only for
real-analytic sub-Riemannian structures while in the general
smooth case we guaranteed the existence of a nonempty open subset
of smooth points. I am indebted to Ludovic Rifford for the nice
observation (see Lemma~1) which allows to get rid of the
analyticity assumption.

\end{document}